\documentstyle[leqno]{article}

\newtheorem{th}{Theorem}[section]
\newtheorem{prop}[th]{Proposition}

\newcounter{defin}[section]
\renewcommand{\thedefin}{\thesection.\arabic{defin}}
\newcounter{ex}[section]

\newcounter{rem}[section]
\renewcommand{\therem}{\thesection.\arabic{rem}}
\title{The Geometry of Autonomous Metrical Multi-Time \\
Lagrange Space of Electrodynamics}
\date{}
\author{Mircea Neagu}
\begin{document}
\maketitle
\begin{abstract}
Section 1 contains physical and geometrical aspects that motivates us to
study the autonomous multi-time Lagrangian space of electrodynamics.
Section 2  constructs the canonical nonlinear connection $\Gamma$ and the
Cartan canonical\linebreak $\Gamma$-linear connection of this space. Section 3
describes the Maxwell equations which govern the electromagnetic
field of this space. The Einstein equations of gravitational potentials of
the autonomous multi-time Lagrange space are written in Section 4. The
conservation laws of these equations will be also described.
\end{abstract}
{\bf Mathematics Subject Classification (2000):} 53B40, 53C60, 53C80.\\
{\bf Key words:} 1-jet fibre bundle, nonlinear connection, Cartan canonical
connection, Maxwell equations, Einstein equations.

\section{Geometrical and physical aspects}

\hspace{5mm} In the last thirty years, many geometrical models in Mechanics or
Physics were based on the notion of ordinary Lagrangian. Thus, the geometrical
concept of Lagrange space was  introduced. The differential geometry of the
Lagrange  spaces is now considerably developped and used in various fields to
study the natural processes where the dependence on position, velocity or
momentum are involved \cite{5}.  We recall that a {\it Lagrange space}
$L^n=(M,L(x,y))$ is defined as a pair which consists of a real, $n$-dimensional
manifold $M$ coordinated by $x=(x^i)_{i=\overline{1,n}}$ and a {\it regular}
Lagrangian $L:TM\to R$ (i. e. the fundamental metrical d-tensor
$\displaystyle{g_{ij}(x,y)={1\over 2}{\partial^2L\over
\partial y^i\partial y^j}}$ is of rank $n$ and has a constante
signature on $TM\backslash\{0\}$). We point  out that  the Lagrangian $L$ is
not necessarily homogenous with respect to the direction $y=(y^i)_{i=\overline{1,n}}$.

An important and well known example of  Lagrange space comes from electrodynamics.
Thus, let us consider the Lagrangian $L:TM\to R$ which governs the movement  law  of a particle
of mass $m\neq  0$ and electric  charge $e$, placed concomitantly into a gravitational
field and an electromagnetic one,
\begin{equation}
L(x,y)=mc\varphi_{ij}(x)y^iy^j+{2e\over m}A_i(x)y^i+U(x),
\end{equation}
where the pseudo-Riemannian metric $\varphi_{ij}(x)$ represents the {\it gravitational
potentials} of the space $M$, $A_i(x)$ are  the components of a covector field
on $M$ representing the {\it  electromagnetic potentials}, $U(x)$ is a function
on $M$  which is called {\it potential function} and $c$ is the physical
constant of light speed. It is obvious that $L$ is a  regular Lagrangian  and,
consequently, the pair $L^n=(M,L)$  is a Lagrange space. This space is known
under the name of {\it the  Lagrange space of electrodynamics}.

At the same time, there are many problems in Physics and Variational Calculus
in which multi-time Lagrangians functions $L$ depending of first order partial
derivatives or, alternatively, of partial directions, are involved. In this
context, the Lagrangian function $L$  is defined on the total space  of the 1-jet
fibre bundle $J^1(T,M)$, where $T$
is a smooth, real, p-dimensional manifold
coordinated by $t=(t^\alpha)_{\alpha=\overline{1,p}}$, whose physical meaning
is that of {\it "multidimensional time"}. We point out that $J^1(T,M)$ is
coordinated by $(t^\alpha,x^i,x^i_\alpha)$.

It is well known that the jet fibre
bundle of order one  is a basic object in the study of classical and quantum
field theories. From a certain physical point of view, the 1-jet fibre bundle
$J^1(T,M)\to T\times M$ is regarded as a {\it bundle of configurations}
\cite{3}, \cite{10}, \cite{11}.

Let us consider the particular case $T=R$ (i. e. the  usual time
axis represented by the set of real  numbers) in the  construction of jet bundle
of order one, we find the bundle
\begin{equation}
J^1(R,M)\equiv R\times TM\to R\times M,\quad (t,x^i,y^i)\to (t,x^i),
\end{equation}
that is, the bundle of configurations of {\it
relativistic rheonomic mechanics} whose invariance gauge group is \cite{10}
\begin{equation}\label{G_1}
\left\{\begin{array}{l}
\tilde t=\tilde t(t)\\
\tilde x^i=\tilde x^i(x^j)\\
\displaystyle{\tilde y^i={\partial\tilde x^i\over\partial x^j}{dt\over
d\tilde t}y^j},
\end{array}\right.
\end{equation}
where the coordinates $(t^1,x^i,x^i_1)$ of the jet space $J^1(R,M)\equiv R
\times TM$ are denoted by $(t,x^i,y^i)$. It is obvious that the form of this
gauge group emphasizes the {\it relativistic} character of the time coordinate
$t$.

We underline that, in the {\it classical rheonomic mechanics} studied in
\cite{5}, the bundle of configurations is
\begin{equation}
R\times TM\to M,\quad (t,x^i,y^i)\to (x^i),
\end{equation}
whose geometrical invariance group is of the form
\begin{equation}\label{G_2}
\left\{\begin{array}{l}
\tilde t=t\\
\tilde x^i=\tilde x^i(x^j)\\
\displaystyle{\tilde y^i={\partial\tilde x^i\over\partial x^j}y^j,}
\end{array}\right.
\end{equation}
that is, it ignores the temporal reparametrizations. Consequently, in  that case,
the temporal coordinate $t$ has a character of {\it absolute} time.

We emphasize that, in the relativistic rheonomic mechanics, a basic role is played by the following
Lagrangian of {\it relativistic rheonomic electrodynamics},
\begin{equation}
{\cal L}=\left[mc\psi^{11}(t)\varphi_{ij}(x)y^iy^j+{2e\over m}A^{(1)}_{(i)}
(t,x)y^i+U(t,x)\right]\sqrt{\vert \psi_{11}\vert},
\end{equation}
where $\psi_{11}$ is a semi-Riemannian  metric on $R$, $A^{(1)}_{(i)}(t,x)$ is
a distiguished tensor on $J^1(R,M)$ and $U(t,x)$ is a smooth function on $R\times M$.

At the same time, we point out that the Lagrangian which governs the {\it classical
rheonomic electrodynamics} has the form
\begin{equation}
L=mc\varphi_{ij}(x)y^iy^j+{2e\over m}A_i(t,x)y^i+U(t,x).
\end{equation}

It is important to note the difference between the notions of Lagrangian
used in both relativistic and classical rheonomic mechanics. From this point of
view, the reader is invited to compare them, following the expositions done in
\cite{5}, \cite{10}. Thus, according to Olver terminology \cite{4}, we point out
that, in the background of relativistic rheonomic mechanics, a  Lagrangian
${\cal L}$ on $J^1(R,M)$ is a  local function on the 1-jet space, which transforms
by the rule $\tilde{\cal L}={\cal L}\vert dt/d\tilde t\vert$. The notion of
Lagrangian function $L$ \linebreak
(i. e. , a smooth function $L:J^1(R,M)\to R$) is also
involved in relativistic rheonomic mechanics. We point out that, in that case,
the geometrical invariance group of the bundle of configurations $J^1(R,M)\to
R\times M$ is \ref{G_1}. In contrast, in the classical rheonomic mechanics, a
Lagrangian $L$ is only a smooth function on the total space of the bundle $R\times TM
\to M$. We remark that, in that case, the bundle of configurations has the
geometrical invariance group \ref{G_2}.

Now, returning us to the general multi-temporal context, we point out that a
fundamental geometrical concept used  in the geometrization of a multi-time
Lagrangian  is that of {metrical multi-time Lagrange space}, introduced in \cite{12}.
The differential geometry of metrical multi-time Lagrange spaces is now considerably
developped in \cite{8}, \cite{12}.

In order to develope this geometrical approach, we fix a semi-Riemannian metric
$\psi=\psi_{\alpha\beta}(t^\gamma)$ on the temporal manifold $T$. We recall that a {\it metrical
multi-time Lagrange space} is a pair $ML^n_p=(J^1(T,M),L)$ consisting of  1-jet
space and a {\it Kronecker $\psi$-regular} multi-time Lagrange function $L$, that
is \cite{12}
\begin{equation}
G^{(\alpha)(\beta)}_{(i)(j)}(t^\gamma,x^k,x^k_\gamma)={1\over 2}{\partial^2L
\over\partial x^i_\alpha\partial x^j_\beta}=\psi^{\alpha\beta}(t^\gamma)
\varphi_{ij}(t^\gamma,x^k,x^k_\gamma),
\end{equation}
where $\varphi_{ij}(t^\gamma,x^k,x^k_\gamma)$ is a d-tensor on $J^1(T,M)$,
symmetric, of rank $n$ and having a constant signature.

An important example of metrical multi-time Lagrange space, which comes from
physics, is offered by the {\it "energy"} Lagrangian function $L$ used in the
Polyakov model of bosonic strings,
\begin{equation}
L(t^\gamma,x^k,x^k_\gamma)={1\over 2}\psi^{\alpha\beta}(t)\varphi_{ij}(x)x^i_
\alpha x^j_\beta.
\end{equation}
We recall that the extremals of the Lagrangian ${\cal L}=L\sqrt{\vert\psi
\vert}$  are exactly the harmonic maps between the semi-Riemannian spaces
$(T,\psi)$ and $(M,\varphi)$.

By a  natural extension of previous examples  of Lagrangian functions, we can
offer another important example of metrical multi-time Lagrange space, considering
the general Lagrangian function $L$ which comes from electrodynamics and theory of
bosonic strings, namely,
\begin{equation}
L=mc\psi^{\alpha\beta}(t)\varphi_{ij}(x)x^i_\alpha x^j_\beta+
{2e\over m}A^{(\alpha)}_{(i)}(t,x)x^i_\alpha+U(t,x),
\end{equation}
where $A^{(\alpha)}_{(i)}(t,x)$ is a distiguished tensor on $J^1(T,M)$ and
$U(t,x)$ is a smooth function on $T\times M$.

In this context,  in order to unify all Lagrangian entities exposed above,
we introduce the following\medskip\\
\addtocounter{defin}{1}
{\bf Definition \thedefin} The pair $EDML^n_p=(J^1(T,M),L)$ which consists of
jet fibre bundle of order one and a Lagrangian function of the form
\begin{equation}
L(t^\gamma,x^k,x^k_\gamma)=h^{\alpha\beta}(t^\gamma)g_{ij}(x^k)x^i_\alpha
x^j_\beta+U^{(\alpha)}_{(i)}(t^\gamma,x^k)x^i_\alpha+F(t^\gamma,x^k)
\end{equation}
where $h_{\alpha\beta}(t^\gamma)$ (resp. $g_{ij}(x^k)$) is a semi-Riemannian
metric on  the temporal (resp. spatial) manifold $T$ (resp. $M$), $U^{(i)}_
{(\alpha)}(t^\gamma,x^k)$ are the local components of  a distinguished tensor
on $J^1(T,M)$ and $F(t^\gamma,x^k)$ is a smooth function on the product manifold
$T\times M$, is called an {\it autonomous metrical multi-time Lagrange space
of electrodynamics}.\medskip\\
\addtocounter{rem}{1}
{\bf Remark \therem} We point out that the non-dynamical character (i. e. ,
the independence with respect to the temporal coordinates) of the spatial metric $g_{ij}
(x^k)$ determined us to use the terminology of {\it autonomous} in the previous
definition.

\section{The geometry of autonomous metrical multi-time
Lagrange space of electrodynamics}

\setcounter{equation}{0}
\hspace{5mm} In this section, we will apply the general geometrical development
of a metrical multi-time Lagrange space \cite{12}, to the particular space of
electrodynamics $EDML^n_p$.

In order to do this development, let us consider the energy action functional associated
to the multi-time Lagrangian of electrodynamics
\begin{equation}
{\cal L}=L\sqrt{\vert h\vert}=\left[h^{\alpha\beta}(t^\gamma)g_{ij}(x^k)x^i_\alpha
x^j_\beta+U^{(\alpha)}_{(i)}(t^\gamma,x^k)x^i_\alpha+F(t^\gamma,x^k)\right]
\sqrt{\vert h\vert},
\end{equation}
namely,
\begin{equation}
{\cal E}_{\cal L}:C^\infty(T,M)\to R,\quad
{\cal E}_{\cal L}(f)=\int_T{\cal L}dt^1\wedge dt^2\ldots\wedge dt^p,
\end{equation}
where the temporal manifold $T$ is considered compact and orientable, the
local expression of the smooth map $f$ is $(t^\alpha)\to(x^i(t^\alpha))$ and
$x^i_\alpha=\partial x^i/\partial t^\alpha$.
In this context, it is proved in \cite{12}
\begin{th}
The extremals of the energy functional  ${\cal E}_L$ associated to the multi-time
Lagrangian
${\cal L}$ are harmonic maps \cite{11} of the spray $(H,G)$ defined by the
temporal components
$$
H^{(i)}_{(\alpha)\beta}=-{1\over 2}H^\gamma_{\alpha\beta}x^i_\gamma
$$
and the local spatial components
$$
G^{(i)}_{(\alpha)\beta}={1\over 2}\gamma^i_{jk}x^j_\alpha x^k_\beta+
{h_{\alpha\beta}g^{il}\over 4p}\left[U^{(\mu)}_{(l)m}x^m_\mu+{\partial U^{(\mu)}
_{(l)}\over\partial t^\mu}+U^{(\mu)}_{(l)}H^\gamma_{\mu\gamma}-{\partial F
\over\partial x^l}\right],
$$
where $H^\gamma_{\alpha\beta}$ (resp. $\gamma^i_{jk}$) are the Christoffel
symbols of the metric $h_{\alpha\beta}$ (resp. $g_{ij}$), $p=\dim T$, and
$\displaystyle{U^{(\alpha)}_{(i)j}={\partial U^{(\alpha)}_{(i)}\over\partial
x^j}-{\partial U^{(\alpha)}_{(j)}\over\partial x^i}}$. In other words, these
extremals verify the harmonic map equations attached to the spray $(H,G)$,
namely,
\begin{equation}
h^{\alpha\beta}\left\{x^i_{\alpha\beta}+2H^{(i)}_{(\alpha)\beta}+2G^{(i)}_
{(\alpha)\beta}\right\}=0.
\end{equation}
\end{th}
\addtocounter{defin}{1}
{\bf Definition \thedefin} The spray $(H,G)$ constructed in  the previous
theorem is called the {\it canonical spray attached to the autonomous metrical
multi-time Lagrange space of electrodynamics}.\medskip

Using the canonical spray $(H,G)$, one naturally induces a nonlinear connection
\linebreak $\Gamma=(M^{(i)}_{(\alpha)\beta}, N^{(i)}_{(\alpha)j})$
on $J^1(T,M)$, which is also called the {\it  canonical nonlinear connection
of the autonomous metrical multi-time Lagrange space of electrodynamics.} Thus,
denoting ${\cal G}^i=h^{\alpha\beta}G^{(i)}_{(\alpha)\beta}$, we establish
the folowing theorem \cite{12}
\begin{th}
The canonical nonlinear connection of the autonomous metrical multi-time Lagrange
space of electrodynamics is determined by the temporal components
\begin{equation}
M^{(i)}_{(\alpha)\beta}=2H^{(i)}_{(\alpha)\beta}=-H^\gamma_{\alpha\beta}x^i_
\gamma
\end{equation}
and the local spatial components
\begin{equation}
N^{(i)}_{(\alpha)j}={\partial{\cal G}^i\over\partial x^j_\gamma}h_{\alpha\gamma}
=\gamma^i_{jk}x^k_\alpha+{h_{\alpha\gamma}g^{il}\over 4}U^{(\gamma)}_{(l)j}.
\end{equation}
\end{th}

Following the paper \cite{12}, by a direct calculation, we determine the {\it
Cartan canonical connection} of the autonomous metrical multi-time Lagrange
space of electrodynamics, together with its torsion and curvature local d-tensors.
\begin{th}
i) The Cartan canonical connection
$
C\Gamma=(H^\gamma_{\alpha\beta}, G^k_{j\gamma}, L^i_{jk}, C^{i(\gamma)}_
{j(k)})
$
of the autonomous metrical multi-time Lagrange space of electrodynamics has
the coefficients
\begin{equation}
H^\gamma_{\alpha\beta}=H^\gamma_{\alpha\beta},\quad G^k_{j\gamma}=0, \quad
L^i_{jk}=\gamma^i_{jk},\quad C^{i(\gamma)}_{j(k)}=0.
\end{equation}

ii) The torsion {\em\bf T} of the Cartan canonical connection of the autonomous
metrical multi-time Lagrange space of electrodynamics is determined by three
local d-tensors, namely,
\begin{equation}
\begin{array}{l}\medskip
\displaystyle{
R^{(m)}_{(\mu)\alpha\beta}=-H^\gamma_{\mu\alpha\beta}x^m_\gamma,
\quad R^{(m)}_{(\mu)\alpha j}=-{h_{\mu\eta}g^{mk}\over 4}\left[H^\eta_
{\alpha\gamma}U^{(\gamma)}_{(k)j}+{\partial U^{(\eta)}_{(k)j}\over\partial
t^\alpha}\right],}\\
\displaystyle{
R^{(m)}_{(\mu)ij}=r^m_{ijk}x^k_\mu+{h_{\mu\eta}g^{mk}\over 4}
\left[U^{(\eta)}_{(k)i\vert j}+U^{(\eta)}_{(k)j\vert i}\right],}
\end{array}
\end{equation}
where $H^\gamma_{\mu\alpha\beta}$ (resp. $r^m_{ijk}$) are the local curvature
tensors of the semi-Riemannian metric $h_{\alpha\beta}$ (resp. $g_{ij}$) and
"$_{\vert i}$" represents the local spatial horizontal covariant derivative
induced by the Cartan connection.

iii) The curvature {\em\bf R} of the Cartan canonical connection of the autonomous
metrical multi-time Lagrange space of electrodynamics is determined by two
local d-tensors, namely, $H^\eta_{\alpha\beta\gamma}$ and $R^l_{ijk}=r^l_{ijk}$,
that is, exactly the curvature tensors of the semi-Riemannian metrics
$h_{\alpha\beta}$ and $g_{ij}$.
\end{th}

\section{Maxwell equations of autonomous metrical multi-time Lagrange space of
electrodynamics}

\setcounter{equation}{0}
\hspace{5mm} In order to develope the electromagnetic theory of the autonomous
metrical multi-time Lagrange space, let us consider the {\it canonical Liouville
d-tensor} {\bf C}=$\displaystyle{x^i_\alpha{\partial\over\partial  x^i_\alpha}}$
on the jet fibre bundle of order one, and let us construct the {\it deflection
d-tensors} \cite{8}
\begin{equation}
\left\{\begin{array}{l}\medskip
\bar D^{(i)}_{(\alpha)\beta}=x^i_{\alpha/\beta}=0\\\medskip
\displaystyle{
D^{(i)}_{(\alpha)j}=x^i_{\alpha\vert j}=-{1\over 4}g^{im}h_{\alpha\mu}U^{(\mu)}
_{(m)j}}\\
d^{(i)(\beta)}_{(\alpha)(j)}=x^i_\alpha\vert^{(\beta)}_{(j)}=\delta^i_j\delta
^\beta_\alpha,
\end{array}\right.
\end{equation}
where "$_{/\beta}$", "$_{\vert j}$" and "$\vert^{(\beta)}_{(j)}$" are the
local covariant derivatives induced by the Cartan canonical connection $C\Gamma$.

Multiplying the deflection d-tensors by the vertical fundamental metrical d-tensor
$h^{\alpha\beta}(t)g_{ik}(x)$, we  construct the d-tensors,
\begin{equation}
\left\{\begin{array}{l}\medskip
\bar D^{(\alpha)}_{(i)\beta}=0\\\medskip
\displaystyle{
D^{(\alpha)}_{(i)j}=-{1\over 4}U^{(\alpha)}_{(i)j}}\\
d^{(\alpha)(\beta)}_{(i)(j)}=h^{\alpha\beta}g_{ij}.
\end{array}\right.
\end{equation}

Taking into  account the general expressions of the local {\it electromagnetic
d-tensors} of a metrical multi-time Lagrange space \cite{8}, by a direct calculation,
we deduce the following
\begin{prop}
The local electromagnetic d-tensors of the autonomous metrical multi-time
Lagrange space of electrodynamics have the expressions,
\begin{equation}
\left\{\begin{array}{l}\medskip
\displaystyle{F^{(\alpha)}_{(i)j}={1\over 2}\left[D^{(\alpha)}_{(i)j}-
D^{(\alpha)}_{(j)i}\right]={1\over 8}\left[U^{(\alpha)}_{(j)i}-U^{(\alpha)}_
{(i)j}\right]=-{1\over 4}U^{(\alpha)}_{(i)j}}\\
\displaystyle{f^{(\alpha)(\beta)}_{(i)(j)}={1\over 2}\left[d^{(\alpha)(\beta)}
_{(i)(j)}-d^{(\alpha)(\beta)}_{(j)(i)}\right]=0.}
\end{array}\right.
\end{equation}
\end{prop}

Particularizing the Maxwell equations of electromagnetic field, described in
the general case of a metrical multi-time Lagrange space \cite{8}, we  deduce
the main result of the electromagnetism of the autonomous metrical multi-time
Lagrange space of electrodynamics.
\begin{th}
The electromagnetic local components $F^{(\alpha)}_{(i)j}$ of the autonomous
metrical multi-time Lagrange space of electrodynamics are governed by the
following equations of Maxwell type,
\begin{equation}
\left\{\begin{array}{l}\medskip
\displaystyle{
F^{(\alpha)}_{(i)j/\beta}={1\over 2}{\cal A}_{\{i,j\}}h^{\alpha\mu}g_{im}
R^{(m)}_{(\mu)\beta j}}\\\medskip
\sum_{\{i,j,k\}}F^{(\alpha)}_{(i)j\vert k}=0\\
\sum_{\{i,j,k\}}F^{(\alpha)}_{(i)j}\vert^{(\gamma)}_{(k)}=0,
\end{array}\right.
\end{equation}
where ${\cal A}_{\{i,j\}}$ represents an alternate sum and $\sum_{\{i,j,k\}}$
means a cyclic sum.
\end{th}

\section{Einstein equations and conservation laws of autonomous metrical
multi-time Lagrange space of electrodynamics}

\setcounter{equation}{0}

\hspace{5mm} In order to develope the gravitational theory of the autonomous
metrical multi-time Lagrange space of electrodynamics $EDML^n_p$, we point out that the
vertical metrical d-tensor $G^{(\alpha)(\beta)}_{(i)(j)}=h^{\alpha\beta}(t)
g_{ij}(x)$ and  the canonical nonlinear connection\linebreak $\Gamma=(M^{(i)}_{(\alpha)
\beta},N^{(i)}_{(\alpha)j})$ of this space induce a natural {\it gravitational
$h$-potential} on the 1-jet space $J^1(T,M)$, which is expressed by \cite{8}
$$
G=h_{\alpha\beta}dt^\alpha\otimes dt^\beta+g_{ij}dx^i\otimes dx^j+h^{\alpha
\beta}g_{ij}\delta x^i_\alpha\otimes\delta x^j_\beta.
$$
Let us consider $C\Gamma=(H^\gamma_{\alpha\beta},G^k_{j\gamma},L^i_{jk},
C^{i(\gamma)}_{j(k)})$ the Cartan canonical connection of $ML^n_p$.

We postulate that the Einstein which govern the gravitational $h$-potential
$G$ of the metrical multi-time Lagrange space of electrodynamics $EDML^n_p$
are the Einstein equations attached to the Cartan canonical connection and
the adapted metric $G$ on $J^1(T,M)$, that is,
$$
Ric(C)-{Sc(C)\over 2}G={\cal K}{\cal T},
$$
where $Ric(C)$ represents the Ricci d-tensor of the Cartan connection, $Sc(C)$
is its scalar curvature, ${\cal K}$ is the Einstein constant and ${\cal T}$
is an intrinsec tensor of matter which is called  the {\it stress-energy}
d-tensor.

In an adapted basis $(X_A)=\displaystyle{\left({\delta\over\delta t^\alpha},
{\delta\over\delta x^i},{\partial\over\partial x^i_\alpha}\right)}$, the curvature
d-tensor {\bf R} of the Cartan connection is expressed locally by {\bf R}$(X_
C,X_B)X_A=R^D_{ABC}X_d$. Hence, it follows that we have $R_{AB}=Ric(X_A,X_B)
=R^D_{ABD}$ and $Sc(C)=G^{AB}R_{AB}$, where
$$
G^{AB}=\left\{\begin{array}{ll}\medskip
h_{\alpha\beta},&\mbox{for}\;\;A=\alpha,\;B=\beta\\\medskip
g^{ij},&\mbox{for}\;\;A=i,\;B=j\\\medskip
h_{\alpha\beta}g^{ij},&\mbox{for}\;\;A={(i)\atop(\alpha)},\;B={(j)\atop(\beta)}\\
0,&\mbox{otherwise}.
\end{array}\right.
$$

Taking into account the expressions of the local curvature d-tensors of the
Cartan connection, we deduce that we have the following two effective local
Ricci d-tensors, namely, $H_{\alpha\beta}$ and $R_{ij}=r_{ij}$, where $H_{\alpha
\beta}$ (resp. $r_{ij}$) are the Ricci tensors associated to the semi-Riemannian
metric $h_{\alpha\beta}$ (resp. $g_{ij}$), the rest of these vanishing.

Denoting $H=h^{\alpha\beta}H_{\alpha\beta}$, $R=g^{ij}R_{ij}$ and $S=h_{\alpha
\beta}g^{ij}R^{(\alpha)(\beta)}_{(i)(j)}$, the scalar curvature of Cartan
connection becomes $Sc(C)=H+R+S$. By a direct calculation, we conclude that the effective scalar curvatures
of a metrical multi-time Lagrange space are $H=h^{\alpha\beta}H_{\alpha\beta}$
and $r=g^{ij}r_{ij}$, that is, exactly the scalar curvatures of the semi-Riemannian
metrics $h_{\alpha\beta}$ and $g_{ij}$.

Consequently, we can establish the following
\begin{th}
The Einstein equations which govern the gravitational $h$-potential $G$, induced
by the Lagrangian function of autonomous metrical multi-time Lagrange space of
electrodynamics, take the form
$$
\left\{\begin{array}{l}\medskip
\displaystyle{H_{\alpha\beta}-{H+r\over 2}h_{\alpha\beta}={\cal K}{\cal T}_
{\alpha\beta}}\\\medskip
\displaystyle{r_{ij}-{H+r\over 2}g_{ij}={\cal K}{\cal T}_{ij}}\\\medskip
\displaystyle{-{H+r\over 2}h^{\alpha\beta}g_{ij}={\cal K}{\cal T}^{(\alpha)
(\beta)}_{(i)(j)}},
\end{array}\right.\leqno{(E_1)}
$$
$$
\left\{\begin{array}{lll}\medskip
0={\cal T}_{\alpha i},&0={\cal T}_{i\alpha},&
0={\cal T}^{(\alpha)}_{(i)\beta}\\
0={\cal T}^{\;(\beta)}_{\alpha(i)},&0={\cal T}^{\;(\alpha)}_{i(j)},&
0={\cal T}^{(\alpha)}_{(i)j}.
\end{array}\right.\leqno{(E_2)}
$$
\end{th}
\addtocounter{rem}{1}
{\bf Remarks \therem} i) Asumming that $p=\dim T>2$ and $n=\dim M>2$, the set
$(E_1)$ of the Einstein equations can be rewritten in the more natural form
$$
\left\{\begin{array}{l}\medskip
\displaystyle{H_{\alpha\beta}-{H\over 2}h_{\alpha\beta}={\cal K}\tilde{\cal T}_
{\alpha\beta}}\\\medskip
\displaystyle{R_{ij}-{R\over 2}g_{ij}={\cal K}\tilde{\cal T}_{ij}},
\end{array}\right.\leqno{(E^\prime_1)}
$$
where $\tilde{\cal T}_{AB},\;A,B\in\{\alpha,i\}$ are the adapted local components
of a new stress-energy d-tensor $\tilde{\cal T}$. This new form of the Einstein
equations will be treated detalied in the more general case of a {\it generalized
metrical multi-time Lagrange space} \cite{9}.

ii) It is remarkable that, writing the Einstein equations of metrical multi-time
Lagrange space of electrodynamics in the new form $(E^\prime_1)$, these reduce to
the classical ones.\medskip

Note that, in order to have the compatibility of the Einstein equations, it is
necessary that the certain  adapted local components of the stress-energy
d-tensor vanish {\it "a priori"}.

At the same time, it is  well known that, from a physical point of view, the
stress-energy d-tensor ${\cal T}$ must  verify the local {\it conservation
laws} ${\cal T}^B_{A\vert B}=0,\;\forall\;A\in\{\alpha,i,{(\alpha)\atop (i)}\}$,
where ${\cal T}^B_A=G^{BD}{\cal T}_{DA}$.

In these conditions, by computations, we obtain the following
\begin{th}
The conservation laws of the Einstein equations of the gravitational $h$-potential
of autonomous metrical multi-time Lagrange space reduce to
$$
\left\{\begin{array}{l}\medskip
\displaystyle{\left[H^\mu_\beta-{H+r\over 2}\delta^\mu_\beta\right]_{/\mu}=0}
\\
\displaystyle{\left[r^m_j-{H+r\over 2}\delta^m_j\right]_{\vert m}=0}.
\end{array}\right.
$$
\end{th}
\addtocounter{rem}{1}
{\bf Remark \therem} Taking into account the components
$\tilde{\cal T}_{\alpha\beta}$ and $\tilde{\cal T}_{ij}$ of the new stress-energy
d-tensor ${\cal T}$ appeared in the $(E_1^\prime)$ form of the Einstein
equations, the conservation laws modify in the following simple and natural
form
\begin{equation}
\tilde{\cal T}^\mu_{\beta/\mu}=0,\;\tilde{\cal T}^m_{j\vert m}=0.
\end{equation}
All these will be discussed detalied in \cite{9}.

\begin{center}
University POLITEHNICA of Bucharest\\
Department of Mathematics I\\
Splaiul Independentei 313\\
77206 Bucharest, Romania\\
e-mail: mircea@mathem.pub.ro\\
\end{center}

\end{document}